\documentclass[12pt]{amsart}
\def\N{{\Bbb N}}
\def\Z{{\Bbb Z}}

\newtheorem{Theorem}{Theorem}[section]
\newtheorem{Corollary}[Theorem]{Corollary}
\newtheorem{Lemma}[Theorem]{Lemma}
\newtheorem{Proposition}[Theorem]{Proposition}
\newtheorem{Thm}{Theorem}
\newtheorem{Cor}[Thm]{Corollary}
\theoremstyle{definition}
\newtheorem{Definition}[Theorem]{Definition}

\newtheorem{Question}[Theorem]{Question}
\theoremstyle{remark}
\newtheorem*{Remark}{Remark}
\begin{document}
\sloppy
\title{Minimality of the boundary of a right-angled Coxeter system}
\author{Tetsuya Hosaka} 
\address{Department of Mathematics, Faculty of Education, 
Utsunomiya University, Utsunomiya, 321-8505, Japan}
\date{November 11, 2006}
\email{hosaka@cc.utsunomiya-u.ac.jp}
\thanks{
Partly supported by the Grant-in-Aid for Young Scientists (B), 
The Ministry of Education, Culture, Sports, Science and Technology, Japan.
(No.\ 18740025).}
\maketitle
\begin{abstract}
In this paper, 
we show that the boundary $\partial\Sigma(W,S)$ of 
a right-angled Coxeter system $(W,S)$ is minimal 
if and only if $W_{\tilde{S}}$ is irreducible, 
where $W_{\tilde{S}}$ is the minimum parabolic subgroup 
of finite index in $W$.
We also provide several applications and remarks.
In particular, 
we obtain that for a right-angled Coxeter system $(W,S)$, 
the set $\{w^{\infty}\,|\,w\in W, o(w)=\infty\}$ is dense 
in the boundary $\partial\Sigma(W,S)$.
\end{abstract}

%%%%%%%%%%%%%
% Section 1 %
%%%%%%%%%%%%%
\section{Introduction and preliminaries}

The purpose of this paper is to study 
dense subsets of the boundary of a Coxeter system.
A {\it Coxeter group} is a group $W$ having a presentation
$$\langle \,S \, | \, (st)^{m(s,t)}=1 \ \text{for}\ s,t \in S \,
\rangle,$$ 
where $S$ is a finite set and 
$m:S \times S \rightarrow \N \cup \{\infty\}$ is a function 
satisfying the following conditions:
\begin{enumerate}
\item[(1)] $m(s,t)=m(t,s)$ for each $s,t \in S$,
\item[(2)] $m(s,s)=1$ for each $s \in S$, and
\item[(3)] $m(s,t) \ge 2$ for each $s,t \in S$ such that $s\neq t$.
\end{enumerate}
The pair $(W,S)$ is called a {\it Coxeter system}.
If, in addition, 
\begin{enumerate}
\item[(4)] $m(s,t) = 2$ or $\infty$ for each $s,t \in S$ such that $s\neq t$,
\end{enumerate}
then $(W,S)$ is said to be {\it right-angled}.
Let $(W,S)$ be a Coxeter system.
Then $W$ has the {\it word metric} $d_{\ell}$ 
defined by $d_{\ell}(w,w')=\ell(w^{-1}w')$ for each $w,w'\in W$, 
where $\ell(w)$ is the word length of $w$ with respect to $S$.
For a subset $T \subset S$, 
$W_T$ is defined as the subgroup of $W$ generated by $T$, 
and called a {\it parabolic subgroup}.
If $T$ is the empty set, then $W_T$ is the trivial group.
A subset $T\subset S$ is called a {\it spherical subset of $S$}, 
if the parabolic subgroup $W_T$ is finite.

Every Coxeter system $(W,S)$ determines 
a {\it Davis complex} $\Sigma(W,S)$ 
which is a CAT(0) geodesic space (\cite{D1}, \cite{D2}, \cite{D3}, \cite{M}).
Here the $1$-skeleton of $\Sigma(W,S)$ is 
the Cayley graph of $W$ with respect to $S$.
The natural action of $W$ on $\Sigma(W,S)$ is proper, cocompact and by isometries.
If $W$ is infinite, then $\Sigma(W,S)$ is noncompact and 
$\Sigma(W,S)$ can be compactified by adding its ideal boundary
$\partial \Sigma(W,S)$ (\cite{BH}, \cite[\S 4]{D2}).
This boundary 
$\partial \Sigma(W,S)$ is called the {\it boundary of} $(W,S)$.
We note that the natural action of $W$ on $\Sigma(W,S)$ 
induces an action of $W$ on $\partial \Sigma(W,S)$ by homeomorphisms.

A subset $A$ of a space $X$ is said to be {\it dense} in $X$, 
if $\overline{A}=X$. 
A subset $A$ of a metric space $X$ is said to be {\it quasi-dense}, 
if there exists $N>0$ such that 
each point of $X$ is $N$-close to some point of $A$.
Suppose that 
a group $G$ acts on a compact metric space $X$ by homeomorphisms.
Then $X$ is said to be {\it minimal}, 
if every orbit $Gx$ is dense in $X$.

For a negatively curved group $\Gamma$ and 
the boundary $\partial\Gamma$ of $\Gamma$, 
we know that each orbit $\Gamma\alpha$ is dense in $\partial\Gamma$ 
for any $\alpha\in\partial\Gamma$, 
that is, $\partial\Gamma$ is minimal (\cite{Gr}).
We note that Coxeter groups are non-positive curved groups and 
not negatively curved groups in general.
Indeed, there exist examples of Coxeter systems 
whose boundaris are not minimal (cf.\ \cite{Ho}, \cite{Ho3}).
The purpose of this paper is to investigate 
when the boundary of a Coxeter system is minimal.

In \cite[Theorem~1]{Ho}, 
we obtained a sufficient condition of a Coxeter system $(W,S)$ 
such that {\it some} orbit of the Coxeter group $W$ is dense 
in the boundary $\partial\Sigma(W,S)$.
After some preliminaries in Section~2, 
we first show that 
the boundary of such Coxeter system is minimal, 
that is, we prove the following theorem in Section~3.

\begin{Thm}\label{Thm1}
Let $(W,S)$ be a Coxeter system.
Suppose that 
$W^{\{s_0\}}$ is quasi-dense in $W$ 
with respect to the word metric and 
$o(s_0t_0)=\infty$ for some $s_0,t_0\in S$, 
where $o(s_0t_0)$ is the order of $s_0t_0$ in $W$.
Then 
\begin{enumerate}
\item[(1)] $\partial\Sigma(W,S)$ is minimal, and
\item[(2)] $\{w^\infty\,|\, w\in W,\ o(w)=\infty\}$ 
is dense in $\partial\Sigma(W,S)$.
\end{enumerate}
\end{Thm}

Here 
$W^{\{s_0\}}=\{w\in W\,|\, \ell(wt)>\ell(w) \ 
\text{for each}\ t\in S\setminus\{s_0\}\}\setminus\{1\}$.

In Sections~4 and 5, 
we investigate right-angled Coxeter groups 
and we prove the following main theorem.

\begin{Thm}\label{Thm}
For a right-angled Coxeter system $(W,S)$, 
the boundary $\partial\Sigma(W,S)$ is minimal 
if and only if $W_{\tilde{S}}$ is irreducible.
\end{Thm}

Here $W_{\tilde{S}}$ is the minimum parabolic subgroup 
of finite index in $(W,S)$, that is, 
for the irreducible decomposition 
$W=W_{S_1}\times\dots\times W_{S_n}$, 
$\tilde{S}=\bigcup\{S_i\,|\,\text{$W_{S_i}$ is infinite}\}$ (\cite{Ho0}).

We provide several applications of Theorem~\ref{Thm} in Sections~5 and 6.
In particular, we obtain the following corollary.

\begin{Cor}
For a right-angled Coxeter system $(W,S)$, 
the set $\{w^\infty\,|\,w\in W,\ o(w)=\infty\}$ 
is dense in the boundary $\partial\Sigma(W,S)$.
\end{Cor}

In Section~6, 
we provide some remarks on dense subsets of boundaries of CAT(0) groups.

%%%%%%%%%%%%%
% Section 2 %
%%%%%%%%%%%%%
\section{Lemmas on Coxeter groups}

In this section, 
we prove some lemmas for (right-angled) Coxeter groups which are used later.

We first provide some definitions.

\begin{Definition}
Let $(W,S)$ be a Coxeter system and $w\in W$.
A representation $w=s_1\cdots s_l$ ($s_i \in S$) is said to be 
{\it reduced}, if $\ell(w)=l$, 
where $\ell(w)$ is the minimum length of 
word in $S$ which represents $w$.
\end{Definition}

\begin{Definition}
Let $(W,S)$ be a Coxeter system.
For each $w \in W$, 
we define $S(w)= \{s \in S \,|\, \ell(ws) < \ell(w)\}$.
For a subset $T \subset S$, 
we also define $W^T= \{w \in W \,|\, S(w)=T \}$. 
\end{Definition}

The following lemma is known.

\begin{Lemma}[\cite{Bo}, \cite{Br}, \cite{Hu}]\label{lem0}
Let $(W,S)$ be a Coxeter system.
\begin{enumerate}
\item[(1)] Let $w\in W$ and 
let $w=s_1\cdots s_l$ be a representation.
If $\ell(w)<l$, then $w=s_1\cdots \hat{s_i} \cdots \hat{s_j} \cdots s_l$ 
for some $1\le i<j\le l$. 
\item[(2)] For each $w \in W$ and $s \in S$, 
$\ell(ws)$ equals either $\ell(w)+1$ or $\ell(w)-1$, and
$\ell(sw)$ also equals either $\ell(w)+1$ or $\ell(w)-1$.
\item[(3)] For each $w\in W$, 
$S(w)$ is a spherical subset of $S$, i.e., $W_{S(w)}$ is finite.
\end{enumerate}
\end{Lemma}

We can obtain the following lemma from the proof of \cite[Lemma~2.5]{Ho}.

\begin{Lemma}\label{lem1}
Let $(W,S)$ be a Coxeter system and 
let $T$ be a maximal spherical subset of $S$.
Then $W^T$ is quasi-dense in $W$.
\end{Lemma}

\begin{proof}
Let $w\in W$.
There exists an element $w'$ of longest length in the coset $wW_T$.
Then $S(w')=T$ by the proof of \cite[Lemma~2.5]{Ho}.
Here $d_{\ell}(w,w')\le \max\{\ell(v)\,|\,v\in W_T\}$.
Thus $W^T$ is quasi-dense in $W$.
\end{proof}

\begin{Lemma}[{\cite[Lemma~2.3~(3)]{Ho}}]\label{lem11}
Let $(W,S)$ be a Coxeter system and 
$s,t\in S$ such that $o(st)=\infty$.
Then $W^{\{s\}}t\subset W^{\{t\}}$.
\end{Lemma}

We provide some lemmas for right-angled Coxeter groups.
We note that right-angled Coxeter groups are {\it rigid}, 
that is, 
a right-angled Coxeter group determines its Coxeter system uniquely 
up to isomorphism (\cite{R}).

By a consequence of 
Tits' solution to the word problem (\cite{T}, \cite[p.50]{Br}), 
we can obtain the following lemma (cf.\ \cite[Lemma~5]{Ho00}).

\begin{Lemma}\label{lem00}
Let $(W,S)$ be a right-angled Coxeter system, 
let $w\in W$, 
let $w=s_1\cdots s_l$ be a reduced representation 
and let $t,t'\in S$. 
If $tw=t(s_1\cdots s_l)$ is reduced and $tw t'=w$, then 
$t=t'$ and $t s_i=s_i t$ for each $i\in\{1,\dots,l\}$.
\end{Lemma}

Using Lemma~\ref{lem00}, we prove the following lemma.

\begin{Lemma}\label{lem}
Let $(W,S)$ be a right-angled Coxeter system, 
let $U$ be a spherical subset of $S$, 
let $s_0\in S\setminus U$ and 
let $T=\{t\in U\,|\,o(s_0t)=2\}$.
Then $W^Us_0\subset W^{T\cup\{s_0\}}$.
\end{Lemma}

\begin{proof}
Let $w\in W^U$.
To prove that $ws_0\in W^{T\cup\{s_0\}}$, 
we show that $S(ws_0)=T\cup\{s_0\}$.
We note that $\ell(ws_0)=\ell(w)+1$ since $s_0\not\in U=S(w)$.
Hence $s_0\in S(ws_0)$.
Also for each $t\in T$, by the definition of $T$, 
$\ell(ws_0t)=\ell(wts_0)<\ell(ws_0)$, and $t\in S(ws_0)$.
Thus $T\cup\{s_0\}\subset S(ws_0)$.
Next we show that $S(ws_0)\subset T\cup\{s_0\}$.
Let $t\in S(ws_0)$.
Then $\ell(ws_0t)<\ell(ws_0)$.
If $w=a_1\dots a_l$ is a reduced representation, 
then by Lemma~\ref{lem0}~(1), 
$$ ws_0t=(a_1\dots a_l)s_0t=(a_1\dots \hat{a_i}\dots a_l)s_0 $$
for some $i\in\{1,\dots,l\}$, or $t=s_0$.
By Lemma~\ref{lem00}, we obtain that $s_0t=ts_0$.
This implies that if $t\neq s_0$ then $\ell(wt)<\ell(w)$, i.e., $t\in S(w)=U$.
Since $t\in U$ and $s_0t=ts_0$, $t\in T$.
Hence $S(ws_0)\subset T\cup\{s_0\}$.
Thus $S(ws_0)=T\cup\{s_0\}$ and $ws_0\in W^{T\cup\{s_0\}}$.
We obtain that $W^Us_0\subset W^{T\cup\{s_0\}}$.
\end{proof}

The following lemma is known.

\begin{Lemma}[\cite{Bo}, \cite{Hu}]\label{lem2}
For a right-angled Coxeter system $(W,S)$, 
the following statements are equivalent.
\begin{enumerate}
\item[(1)] $(W,S)$ is irreducible.
\item[(2)] For each $a,b\in S$ such that $a\neq b$, 
there exists a sequence $\{a=s_1,s_2,\dots,s_n=b\}\subset S$ 
such that $o(s_is_{i+1})=\infty$ for any $i\in\{1,\dots,n-1\}$.
\end{enumerate}
\end{Lemma}

%%%%%%%%%%%%%
% Section 3 %
%%%%%%%%%%%%%
\section{Minimality of the boundary of a Coxeter system}

In this section, 
we provide an extension of a result in \cite{Ho} 
on minimality of the boundary of a Coxeter system.

\begin{Theorem}\label{thm1}
Let $(W,S)$ be a Coxeter system.
Suppose that 
$W^{\{s_0\}}$ is quasi-dense in $W$ and 
$o(s_0t_0)=\infty$ for some $s_0,t_0\in S$.
Then 
\begin{enumerate}
\item[(1)] $\partial\Sigma(W,S)$ is minimal, and
\item[(2)] $\{w^\infty\,|\, w\in W,\ o(w)=\infty\}$ 
is dense in $\partial\Sigma(W,S)$.
\end{enumerate}
\end{Theorem}

\begin{proof}
Suppose that 
$W^{\{s_0\}}$ is quasi-dense in $W$ and 
$o(s_0t_0)=\infty$ for some $s_0,t_0\in S$.
Then we show that 
$W\alpha$ is dense in $\partial\Sigma(W,S)$ 
for any $\alpha\in\partial\Sigma(W,S)$.

Let $\alpha\in\partial\Sigma(W,S)$ and 
let $\{w_i\}\subset W$ be a sequence which converges to $\alpha$ 
in $\Sigma(W,S)\cup\partial\Sigma(W,S)$.
Since $W^{\{s_0\}}$ is quasi-dense in $W$, 
there exists a number $N>0$ such that 
for each $w\in W$, $d_{\ell}(w,v)\le N$ for some $v\in W^{\{s_0\}}$.
Hence for each $w\in W$, there exists $x\in W$ such that 
$\ell(x)\le N$ and $wx\in W^{\{s_0\}}$.
For each $i$, 
there exists $x_i\in W$ such that 
$\ell(x_i)\le N$ and $(w_i)^{-1}x_i\in W^{\{s_0\}}$.
We note that the set $\{x\in W\,|\, \ell(x)\le N\}$ is finite 
because $S$ is finite.
Hence $\{x_i\,|\,i\in\N\}$ is finite, and 
there exist $x\in W$ and a sequence $\{i_j\,|\,j\in\N\}\subset\N$ 
such that $x_{i_j}=x$ for each $j\in\N$.
Then for each $j\in\N$, $(w_{i_j})^{-1}x_{i_j}=(w_{i_j})^{-1}x \in W^{\{s_0\}}$ and 
$(w_{i_j})^{-1}xt_0 \in W^{\{t_0\}}$ by Lemma~\ref{lem11}, 
since $o(s_0t_0)=\infty$.
Hence $t_0x^{-1}w_{i_j}\in (W^{\{t_0\}})^{-1}$.
The sequence $\{t_0x^{-1}w_{i_j}\,|\,j\in\N\}$ converges to $t_0x^{-1}\alpha$, 
since $\{w_{i_j}\,|\,j\in\N\}$ converges to $\alpha$.
By the proof of \cite[Theorem~4.1]{Ho}, 
we obtain that $Wt_0x^{-1}\alpha$ is dence in $\partial\Sigma(W,S)$, 
that is, $W\alpha$ is dence in $\partial\Sigma(W,S)$.

Thus every orbit $W\alpha$ is dense in $\partial\Sigma(W,S)$ 
and $\partial\Sigma(W,S)$ is minimal.

The minimality of $\partial\Sigma(W,S)$ implies that 
the set $\{w^\infty\,|\,w\in W,\ o(w)=\infty\}$ 
is dense in $\partial\Sigma(W,S)$ (see Proposition~\ref{prop}).
\end{proof}

Here we have a question whether conversely 
if $\partial\Sigma(W,S)$ is minimal then 
$W^{\{s_0\}}$ is quasi-dense in $W$ and 
$o(s_0t_0)=\infty$ for some $s_0,t_0\in S$.
The answer of this question is no in general.

For example, 
let $S=\{s_1,s_2,s_3\}$ and let 
$$W=\langle S\,|\,s_1^2=s_2^2=s_3^2=(s_1s_2)^4=(s_2s_3)^4=(s_3s_1)^4=1\rangle.$$
Then $W$ is a negatively curved group and 
the boundary $\partial\Sigma(W,S)$ is minimal.
On the other hand, 
there do not exist $s_0,t_0\in S$ such that $o(s_0t_0)=\infty$.

In Section~5, we will show that 
the answer of this question is yes for right-angled Coxeter groups.

%%%%%%%%%%%%%
% Section 4 %
%%%%%%%%%%%%%
\section{Key Lemma}

In this section, 
we prove the following lemma which plays a key role in the proof of the main theorem.

\begin{Lemma}\label{KeyLem}
Let $(W,S)$ be a right-angled Coxeter system such that $W$ is infinite.
If $W$ is irreducible, then 
$W^{\{s_0\}}$ is quasi-dense in $W$ for some $s_0\in S$.
\end{Lemma}

\begin{proof}
We suppose that $W^{\{s\}}$ is not quasi-dense in $W$ for any $s\in S$.
Then we show that $W$ is not irreducible.

Let $s_0\in S$, let $T_1=\{t\in S\,|\, o(s_0t)=2\}$ 
and let $S_1=S\setminus T_1$.
If $T_1=\emptyset$ then $o(s_0s)=\infty$ for each $s\in S\setminus\{s_0\}$, 
hence $W^{\{s_0\}}$ is quasi-dense in $W$ which contradicts the assumption.
Thus $T_1\neq\emptyset$.
If $S_1=\{s_0\}$ then $W=W_{\{s_0\}}\times W_{T_1}$, i.e., 
$W$ is not irreducible.
We suppose that $S_1\neq\{s_0\}$.

Let $s_1\in S_1\setminus\{s_0\}$, 
let $T_2=\{t\in T_1\,|\, o(s_1t)=2\}$ 
and let $S_2=S\setminus T_2=S_1\cup(T_1\setminus T_2)$.
We note that $o(s_it)=2$ for each $i\in\{0,1\}$ and $t\in T_2$, 
i.e., $W_{\{s_0,s_1\}\cup T_2}=W_{\{s_0,s_1\}}\times W_{T_2}$.
Since $s_1\in S_1\setminus\{s_0\}$, 
we obtain that $o(s_0s_1)=\infty$ and $W_{\{s_0,s_1\}}$ is irreducible.

Now we show that $T_2\neq\emptyset$.
Suppose that $T_2=\emptyset$.
This means that $o(s_1t)=\infty$ for each $t\in T_1$.
Let $U$ be a maximal spherical subset of $S$ such that $s_0\in U$.
Then $o(uv)=2$ for each $u,v\in U$ such that $u\neq v$, 
because $(W,S)$ is right-angled and $W_U$ is finite.
Hence $o(s_0u)=2$ for each $u\in U$, since $s_0\in U$.
This means that $U\subset T_1\cup\{s_0\}$.
Hence $o(s_1u)=\infty$ for any $u\in U$, 
because $o(s_1t)=\infty$ for any $t\in T_1$ and $o(s_0s_1)=\infty$.
Thus $W^Us_1\subset W^{\{s_1\}}$ by Lemma~\ref{lem}.
Here by Lemma~\ref{lem1}, 
$W^U$ is quasi-dense in $W$, since $U$ is a maximal spherical subset of $S$.
Hence $W^{\{s_1\}}$ is quasi-dense in $W$.
This contradicts the assumption.
Thus we obtain that $T_2\neq\emptyset$.

If $S_2=\{s_0,s_1\}$ then $W=W_{\{s_0,s_1\}}\times W_{T_2}$ and 
$W$ is not irreducible.
We suppose that $S_2\neq\{s_0,s_1\}$.
Let $s_2\in S_2\setminus\{s_0,s_1\}$, 
let $T_3=\{t\in T_2\,|\, o(s_2t)=2\}$ 
and let $S_3=S\setminus T_3=S_2\cup(T_2\setminus T_3)$.

By induction, we define $s_k,T_{k+1},S_{k+1}$ as follows:
Let 
\begin{align*}
&s_k\in S_k\setminus\{s_0,\dots,s_{k-1}\}, \\
&T_{k+1}=\{t\in T_k\,|\,o(s_kt)=2\} \ \text{and} \\
&S_{k+1}=S\setminus T_{k+1}.
\end{align*}
Then $W_{\{s_0,s_1,\dots,s_k\}\cup T_{k+1}}=
W_{\{s_0,s_1,\dots,s_k\}}\times W_{T_{k+1}}$.
If $S_{k+1}\setminus\{s_0,s_1,\dots,s_k\}=\emptyset$ 
then $W=W_{S_{k+1}}\times W_{T_{k+1}}$, i.e., $W$ is not irreducible.
Here we note that $T_{k+1}\subset T_k\subset\dots\subset T_2\subset T_1$.
If $T_{k}\neq\emptyset$ for each $k$, 
then by the finiteness of $S$, 
there exists a number $n$ such that $W=W_{S_n}\times W_{T_n}$, 
hence $W$ is not irreducible.

We prove the following statements by the induction on $k$.
\begin{enumerate}
\item[${\rm (i_k)}$] $T_k\neq\emptyset$.
\item[${\rm (ii_k)}$] $W_{\{s_0,\dots,s_{k-1}\}}$ is irreducible.
\item[${\rm (iii_k)}$] There exists a spherical subset $U_k\subset T_k$ 
such that $W^{U_k\cup\{s_i\}}$ is quasi-dense in $W$ for each $i\in\{0,\dots,k-1\}$.
\end{enumerate}

We first consider in the case $k=2$.
The statement ${\rm (i_2)}$ $T_2\neq\emptyset$ was proved in the above.
Also ${\rm (ii_2)}$ holds, since 
$W_{\{s_0,s_{1}\}}=W_{\{s_0\}}*W_{\{s_1\}}$ is irreducible.
We show that the statement ${\rm (iii_2)}$ holds.
Let $U$ be a maximal spherical subset of $S$ such that $s_0\in U$.
Then $W^U$ is quasi-dense in $W$ by Lemma~\ref{lem1}.
Let $U_2=U\cap T_2$.
We note that $U_2=\{t\in U\,|\, o(s_1t)=2\}$.
By Lemma~\ref{lem}, 
$W^Us_1\subset W^{U_2\cup\{s_1\}}$.
Hence $W^{U_2\cup\{s_1\}}$ is quasi-dense in $W$.
(This implies that $U_2\neq\emptyset$ by the assumption.)
Also $W^{U_2\cup\{s_0\}}$ is quasi-dense in $W$, since 
$W^{U_2\cup\{s_1\}}s_0\subset W^{U_2\cup\{s_0\}}$ by Lemma~\ref{lem}.
Thus ${\rm (iii_2)}$ holds.

We suppose that ${\rm (i_k)}$, ${\rm (ii_k)}$ 
and ${\rm (iii_k)}$ hold for some $k\ge 2$.
Then we prove that ${\rm (i_{k+1})}$, ${\rm (ii_{k+1})}$ 
and ${\rm (iii_{k+1})}$ hold.

${\rm (i_{k+1})}$: 
We show that $T_{k+1}\neq\emptyset$.
Suppose that $T_{k+1}=\emptyset$.
If $o(s_ks_i)=2$ for any $i\in\{0,1,\dots,k-1\}$, 
then $s_k\in T_k$ which contradicts the definition of $s_k$.
Hence $o(s_ks_{i_0})=\infty$ for some $i_0\in\{0,1,\dots,k-1\}$.
Since $T_{k+1}=\emptyset$, $o(s_kt)=\infty$ for any $t\in T_k$.
Here $U_k\subset T_k$ and $o(s_kt)=\infty$ for any $t\in U_k$.
Hence $W^{U_k\cup\{s_{i_0}\}}s_k\subset W^{\{s_k\}}$ by Lemma~\ref{lem}.
By ${\rm (iii_k)}$, $W^{U_k\cup\{s_{i_0}\}}$ is quasi-dense in $W$.
Thus $W^{\{s_k\}}$ is also quasi-dense in $W$, 
which contradicts the assumption.
Hence $T_{k+1}\neq\emptyset$.

${\rm (ii_{k+1})}$: 
We show that $W_{\{s_0,\dots,s_{k-1},s_k\}}$ is irreducible.
Now $o(s_ks_{i_0})=\infty$ for some $i_0\in\{0,1,\dots,k-1\}$ 
by the above argument.
Also $W_{\{s_0,\dots,s_{k-1}\}}$ is irreducible by the hypothesis ${\rm (ii_k)}$.
Hence $W_{\{s_0,\dots,s_{k-1},s_k\}}$ is irreducible.

${\rm (iii_{k+1})}$: 
By ${\rm (iii_k)}$, there exists a spherical subset $U_k\subset T_k$ 
such that $W^{U_k\cup\{s_i\}}$ is quasi-dense in $W$ for each $i\in\{0,\dots,k-1\}$.
We define $U_{k+1}=U_k\cap T_{k+1}$, i.e., $U_{k+1}=\{t\in U_k\,|\,o(s_kt)=2\}$.
Here $o(s_ks_{i_0})=\infty$ for some $i_0\in\{0,1,\dots,k-1\}$ 
by the above argument.
Then $W^{U_k\cup\{s_{i_0}\}}s_k\subset W^{U_{k+1}\cup\{s_k\}}$
by Lemma~\ref{lem}.
Hence $W^{U_{k+1}\cup\{s_k\}}$ is quasi-dense in $W$, 
since $W^{U_k\cup\{s_{i_0}\}}$ is so.
Finally we show that 
$W^{U_{k+1}\cup\{s_i\}}$ is quasi-dense in $W$ for each $i\in\{0,\dots,k-1,k\}$.
We note that $W_{\{s_0,\dots,s_{k-1},s_k\}}$ is irreducible by ${\rm (ii_{k+1})}$.
Hence for each $j_0\in\{0,\dots,k-1\}$, 
there exists a sequence 
$\{s_k=a_0,a_1,\dots,a_m=s_{j_0}\}\subset\{s_i\,|\,i=0,1,\dots,k\}$ 
such that $o(a_ia_{i+1})=\infty$ by Lemma~\ref{lem2}.
Then by Lemma~\ref{lem}, 
\begin{align*}
W^{U_{k+1}\cup\{s_k\}}a_1a_2\cdots a_m
&\subset W^{U_{k+1}\cup\{a_1\}}a_2\cdots a_m \\
&\subset\dots\subset W^{U_{k+1}\cup\{a_m\}}=W^{U_{k+1}\cup\{s_{j_0}\}},
\end{align*}
because $o(s_iu)=2$ for each $i\in\{0,1,\dots,k-1,k\}$ and $u\in U_{k+1}$.
Thus $W^{U_{k+1}\cup\{s_{j_0}\}}$ is quasi-dense in $W$.
Hence ${\rm (iii_{k+1})}$ holds.

Thus by the induction on $k$, 
we can define $s_{k-1}, T_k, S_k$ 
which satisfy ${\rm (i_k)}$, ${\rm (ii_k)}$ and ${\rm (iii_k)}$.
Since $S$ is finite, 
there exists a number $n$ such that $S_n=\{s_0,s_1,\dots,s_{n-1}\}$ and 
$W=W_{S_n}\times W_{T_n}$, where $T_n\neq\emptyset$.
Thus $W$ is not irreducible.
\end{proof}

%%%%%%%%%%%%%
% Section 5 %
%%%%%%%%%%%%%
\section{Dense subsets of the boundary of a right-angled Coxeter group}

We obtain the following main theorem 
from Theorem~\ref{thm1} and Lemma~\ref{KeyLem}.

\begin{Theorem}\label{thm}
Let $(W,S)$ be a right-angled Coxeter system such that $W$ is infinite.
Then the following statements are equivalent.
\begin{enumerate}
\item[(1)] $\partial\Sigma(W,S)$ is minimal.
\item[(2)] $W_{\tilde{S}}$ is irreducible.
\item[(3)] $W^{\{s_0\}}$ is quasi-dense in $W$ and $o(s_0t_0)=\infty$ 
for some $s_0,t_0\in S$.
\end{enumerate}
\end{Theorem}

\begin{proof}
$(3)\Rightarrow(1)$: 
If the statement (3) holds, 
then $\partial\Sigma(W,S)$ is minimal by Theorem~\ref{thm1}.

$(1)\Rightarrow(2)$: 
Suppose that $W_{\tilde{S}}$ is not irreducible.
Let $W_{\tilde{S}}=W_{S_1}\times W_{S_2}$, 
where $W_{S_1}$ and $W_{S_2}$ are infinite.
Then $\partial\Sigma(W,S)=\partial\Sigma(W_{\tilde{S}},\tilde{S})$ and 
$\Sigma(W_{\tilde{S}},\tilde{S})
=\Sigma(W_{S_1},S_1)\times \Sigma(W_{S_2},S_2)$.
Here by \cite[Theorem~4.3]{Ho0}, 
$\partial\Sigma(W_{S_1},S_1)$ is $W$-invariant, that is, 
$W\partial\Sigma(W_{S_1},S_1)=\partial\Sigma(W_{S_1},S_1)$.
Thus for $\alpha\in \partial\Sigma(W_{S_1},S_1)$, 
$\overline{W\alpha}\subset \partial\Sigma(W_{S_1},S_1)$.
Hence $\partial\Sigma(W,S)$ is not minimal.
In Section~6, we will provide more general proof (Theorem~\ref{thm6-1}).

$(2)\Rightarrow(3)$: 
Suppose that $W_{\tilde{S}}$ is irreducible.
By Lemma~\ref{KeyLem}, 
$W^{\{s_0\}}\cap W_{\tilde{S}}$ is quasi-dense 
in $W_{\tilde{S}}$ for some $s_0\in\tilde{S}$.
Here $W=W_{\tilde{S}}\times W_{S\setminus\tilde{S}}$ and 
$W_{S\setminus\tilde{S}}$ is finite (see \cite{Ho0}).
Hence $W^{\{s_0\}}$ is quasi-dense in $W$.
Since $W_{\tilde{S}}$ is irreducible, 
$o(s_0t_0)=\infty$ for some $t_0\in\tilde{S}$ by Lemma~\ref{lem2}.
Thus the statement (3) holds.
\end{proof}

There is the following question in \cite{Ho2}.

\begin{Question}\label{Q:para}
Let $(W,S)$ be a Coxeter system.
Is it the case that 
if $(W,S)$ is an irreducible Coxeter system then
$W\partial\Sigma(W_T,T)$ is dense in $\partial\Sigma(W,S)$ 
for any subset $T$ of $S$ such that $W_T$ is infinite? 
\end{Question}

Theorem~\ref{thm} implies that the answer of Question~\ref{Q:para} 
is yes for right-angled Coxeter groups.
Moreover, 
as an application of Theorem~\ref{thm}, 
we obtain the following corollary.

\begin{Corollary}\label{cor5-2}
Let $(W,S)$ be a right-angled Coxeter system and let $T\subset S$.
Then the following statements are equivalent.
\begin{enumerate}
\item[(1)] $W\partial\Sigma(W_T,T)$ is dense in $\partial\Sigma(W,S)$.
\item[(2)] If $W=W_{S_1}\times\dots\times W_{S_n}$ 
is the irreducible decomposition of $W$, 
then $W_{S_i\cap T}$ is infinite 
for each $i\in\{1,\dots,n\}$ such that $W_{S_i}$ is infinite.
\end{enumerate}
\end{Corollary}

\begin{proof}
$(1)\Rightarrow (2)$: 
Let $W=W_{S_1}\times\dots\times W_{S_n}$ 
be the irreducible decomposition of $W$. 
We suppose that there exists $i_0\in\{1,\dots,n\}$ such that 
$W_{S_{i_0}}$ is infinite and $W_{S_{i_0}\cap T}$ is finite.
Let $A_1=S\setminus S_{i_0}$ and $A_2=S_{i_0}$.
Then $W=W_{A_1}\times W_{A_2}$, 
$W_{A_2}$ is infinite and $W_{A_2\cap T}$ is finite.
We note that 
$\partial\Sigma(W_{A_1},A_1)$ is $W$-invariant by \cite[Theorem~4.3]{Ho0}.
Since $W_T=W_{A_1\cap T}\times W_{A_2\cap T}$ and $W_{A_2\cap T}$ is finite, 
$\partial\Sigma(W_T,T)\subset\partial\Sigma(W_{A_1},A_1)$.
Thus 
$$W\partial\Sigma(W_T,T)\subset W\partial\Sigma(W_{A_1},A_1)
=\partial\Sigma(W_{A_1},A_1).$$
Since $W_{A_2}$ is infinite and 
$$\partial\Sigma(W,S)=\partial\Sigma(W_{A_1},A_1)*\partial\Sigma(W_{A_2},A_2), $$
$W\partial\Sigma(W_T,T)$ is not dene in $\partial\Sigma(W,S)$.

$(2)\Rightarrow (1)$: 
Let $W=W_{S_1}\times\dots\times W_{S_n}$ 
be the irreducible decomposition of $W$. 
Suppose that (2) holds.
Then we prove that (1) holds by induction on $n$.

We first consider in the case $n=1$.
Then $W=W_{S_1}$ is irreducible.
Since $W_{S_1\cap T}$ is infinite, $\partial\Sigma(W_T,T)\neq\emptyset$.
Hence 
$W\partial\Sigma(W_T,T)$ is dense in $\partial\Sigma(W,S)$ by Theorem~\ref{thm}.

Next we consider in the case $n>1$.
Let $A_1=S_1\cup\dots\cup S_{n-1}$ and $A_2=S_n$.
Then $W=W_{A_1}\times W_{A_2}$ and $W_T=W_{A_1\cap T}\times W_{A_2\cap T}$.
Here 
\begin{align*}
W\partial\Sigma(W_T,T)
&=W(\partial\Sigma(W_{A_1\cap T},A_1\cap T)*
\partial\Sigma(W_{A_2\cap T},A_2\cap T)) \\
&\supset W_{A_1}\partial\Sigma(W_{A_1\cap T},A_1\cap T)*
W_{A_2}\partial\Sigma(W_{A_2\cap T},A_2\cap T).
\end{align*}
By the inductive hypothesis, 
$W_{A_i}\partial\Sigma(W_{A_i\cap T},A_i\cap T)$ 
is dense in $\partial\Sigma(W_{A_i},A_i)$ 
for each $i=1,2$.
Since 
$$\partial\Sigma(W,S)=\partial\Sigma(W_{A_1},A_1)*\partial\Sigma(W_{A_2},A_2), $$
we obtain that $W\partial\Sigma(W_T,T)$ is dense in $\partial\Sigma(W,S)$.
\end{proof}

Also we can obtain the following corollary from Theorem~\ref{thm}.
The proof in more general case is provided in Section~6.

\begin{Corollary}\label{cor5-1}
For a right-angled Coxeter system $(W,S)$, 
the set $\{w^\infty\,|\,w\in W,\ o(w)=\infty\}$ 
is dense in the boundary $\partial\Sigma(W,S)$.
\end{Corollary}

%%%%%%%%%%%%%
% Section 6 %
%%%%%%%%%%%%%
\section{Remarks on dense subsets of boundaries of CAT(0) groups}

In this section, 
we investigate dense subsets of boundaries of CAT(0) groups.
Definitions and basic properties of CAT(0) spaces and their boundaries 
can be found in \cite{BH}.
A group $\Gamma$ is called a {\it CAT(0) group}, 
if $\Gamma$ acts geometrically (i.e.\ properly and cocompactly by isometries) 
on some CAT(0) space.
For example,
a Coxeter group $W$ acts geometrically on the Davis complex $\Sigma(W,S)$ 
which is a CAT(0) space, and $W$ is a CAT(0) group.

There is the following open problem.

\begin{Question}\label{Q:R}
Suppose that a group $\Gamma$ acts geometrically on a CAT(0) space $X$.
Is it the case that 
the set $\{\gamma^\infty\,|\,\gamma\in\Gamma,\ o(\gamma)=\infty\}$ 
is dense in the boundary $\partial X$?
\end{Question}

We consider relation between this question 
and minimality of boundaries of CAT(0) groups.

First, we note that there is the following proposition.

\begin{Proposition}\label{prop}
Suppose that a group $\Gamma$ acts geometrically on a CAT(0) space $X$.
If there exists $\delta\in\Gamma$ such that 
$o(\delta)=\infty$ and $\Gamma\delta^\infty$ is dense in the boundary $\partial X$, 
then the set $\{\gamma^\infty\,|\,\gamma\in\Gamma,\ o(\gamma)=\infty\}$ 
is dense in $\partial X$.
Hence, 
if the boundary $\partial X$ is minimal, then 
the set $\{\gamma^\infty\,|\,\gamma\in\Gamma,\ o(\gamma)=\infty\}$ 
is dense in $\partial X$.
\end{Proposition}

\begin{proof}
Suppose that $\delta\in\Gamma$ such that 
$o(\delta)=\infty$ and $\Gamma\delta^\infty$ is dense in $\partial X$.
Let $\alpha\in\partial X$.
Since $\Gamma\delta^\infty$ is dense in $\partial X$, 
there exists a sequence $\{\gamma_i\}\subset\Gamma$ 
such that $\{\gamma_i\delta^\infty\}$ converges to $\alpha$ in $\partial X$.
Here for $x_0\in X$ and each $i$, 
the sequence $\{(\gamma_i\delta\gamma_i^{-1})^jx_0\}_j$ 
converges to $\gamma_i\delta^\infty$ in $X\cup \partial X$.
Hence $(\gamma_i\delta\gamma_i^{-1})^\infty=\gamma_i\delta^\infty$ and 
$\{(\gamma_i\delta\gamma_i^{-1})^\infty\}_i$ 
converges to $\alpha$ in $\partial X$.
Thus $\{\gamma^\infty\,|\,\gamma\in\Gamma,\ o(\gamma)=\infty\}$ 
is dense in $\partial X$.

Now we suppose that the boundary $\partial X$ is minimal.
Every CAT(0) group has an element of infinite order (\cite[Theorem~11]{Sw}).
Let $\delta\in\Gamma$ such that $o(\delta)=\infty$.
Then $\Gamma\delta^\infty$ is dense in $\partial X$ 
because $\partial X$ is minimal.
Hence, by the above argument, 
the set $\{\gamma^\infty\,|\,\gamma\in\Gamma,\ o(\gamma)=\infty\}$ 
is dense in $\partial X$.
\end{proof}

We obtain the following proposition 
from some splitting theorems for CAT(0) spaces.

\begin{Proposition}\label{thm:splitting}
Suppose that a group $\Gamma=\Gamma_1\times \Gamma_2$ acts geometrically 
on a CAT(0) space $X$ 
where $\Gamma_1$ and $\Gamma_2$ are infinite.
Then 
$X$ contains a quasi-dense subspace $X'=X_1\times X_2$ 
and there exists a product subgroup 
$\Gamma'_1\times \Gamma'_2$ of finite index in $\Gamma$ 
such that $X_1$ is the convex hull $C(\Gamma'_1x_0)$ for some $x_0\in X$ and 
$\Gamma'_2$ acts geometrically on $X_2$ by projection.
\end{Proposition}

\begin{proof}
By \cite[Lemma~2.1]{Ho000}, 
there exist subgroups $G_1\times A_1$ and $G_2\times A_2$ 
of finite index in $\Gamma_1$ and $\Gamma_2$ respectively 
such that $G_1$ and $G_2$ have finite center and 
$A_i$ is isomorphic to $\Z^{n_i}$ for some $n_i$ ($i=1,2$).

In the case $A_i$ is not trivial for some $i\in\{1,2\}$, 
let $\Gamma'_1=A_i$ and $\Gamma'_2=G_1\times G_2\times A_{3-i}$.
Then by the Flat Torus Theorem (\cite[Theorem~II.7.1]{BH}), 
the proposition holds.

In the case $A_1$ and $A_2$ are trivial, 
let $\Gamma'_1=G_1$ and $\Gamma'_2=G_2$.
Here $G_1$ and $G_2$ have finite center.
By \cite[Theorem~2]{Ho4} and \cite[Corollary~10]{Mo}, 
the proposition holds.
\end{proof}

Concerning non-minimality of boundaries of CAT(0) groups, 
using Proposition~\ref{thm:splitting}, we can prove the following theorem.

\begin{Theorem}\label{thm6-1}
Suppose that a group $\Gamma$ acts geometrically on a CAT(0) space $X$.
If $\Gamma$ contains a subgroup $\Gamma_1\times\Gamma_2$ of finite index 
such that $\Gamma_1$ and $\Gamma_2$ are infinite, 
then the boundary $\partial X$ is not minimal.
\end{Theorem}

\begin{proof}
Let $\Gamma_1\times\Gamma_2$ be a subgroup of finite index in $\Gamma$, 
where $\Gamma_1$ and $\Gamma_2$ are infinite.
Then $\Gamma_1\times\Gamma_2$ acts geometrically on $X$.
By Proposition~\ref{thm:splitting}, 
$X$ contains a quasi-dense subspace $X_1\times X_2$ 
and there exist a product subgroup 
$\Gamma'_1 \times \Gamma'_2$ of finite index in $\Gamma$ 
such that $X_1$ is the convex hull $C(\Gamma'_1x_0)$ for some $x_0\in X$ and 
$\Gamma'_2$ acts geometrically on $X_2$ by projection.

To prove that $\partial X$ is not minimal, 
we show that $\Gamma(\partial X_1)$ is not dense in $\partial X$.

Since $\Gamma'_1\times\Gamma'_2$ is a subgroup of finite index in $\Gamma$, 
there exist a number $n$ and $\{\delta_1,\dots,\delta_n\}\subset\Gamma$ 
such that $\Gamma=\bigcup_{i=1}^{n}\delta_i(\Gamma'_1\times\Gamma'_2)$.

Since $X_1=C(\Gamma'_1x_0)$ is $\Gamma'_1$-invariant, 
$\Gamma'_1(\partial X_1)=\partial X_1$.
For each $\gamma_2\in \Gamma'_2$, 
$\gamma_2 X_1$ and $X_1$ are parallel 
by the proof of splitting theorems (\cite{BH}, \cite{Ho4}, \cite{Mo}), 
hence $\gamma_2(\partial X_1)=\partial X_1$, that is, 
$\Gamma'_2(\partial X_1)=\partial X_1$.
Thus $(\Gamma'_1\times\Gamma'_2)(\partial X_1)=\partial X_1$.

Hence 
\begin{align*}
\Gamma(\partial X_1)&=
(\bigcup_{i=1}^{n}\delta_i(\Gamma'_1\times\Gamma'_2))(\partial X_1)\\
&=\bigcup_{i=1}^{n}(\delta_i(\Gamma'_1\times\Gamma'_2)(\partial X_1))\\
&=\bigcup_{i=1}^{n}(\delta_i(\partial X_1)).
\end{align*}
Here we note that 
$\Gamma(\partial X_1)=\bigcup_{i=1}^{n}(\delta_i(\partial X_1))$ is closed.
Hence 
\begin{align*}
\dim\overline{\Gamma(\partial X_1)}&=\dim\bigcup_{i=1}^{n}(\delta_i(\partial X_1))
=\dim\partial X_1 \\
&<\dim (\partial X_1*\partial X_2)=\dim\partial X.
\end{align*}
Thus $\Gamma(\partial X_1)$ is not dense in $\partial X$.
This implies that $\partial X$ is not minimal.
\end{proof}

Here the author has the following question 
which is the converse of Theorem~\ref{thm6-1}.

\begin{Question}\label{Q1}
Suppose that a group $\Gamma$ acts geometrically on a CAT(0) space $X$.
Is it the case that if 
$\Gamma$ does not contain a subgroup $\Gamma_1\times\Gamma_2$ of finite index 
such that $\Gamma_1$ and $\Gamma_2$ are infinite, 
then the boundary $\partial X$ is minimal?
\end{Question}

Theorem~\ref{thm} implies that 
the answer of Question~\ref{Q1} is yes 
for right-angled Coxeter groups and their boundaries.

\begin{Remark}
If $\Gamma$ does not contain a subgroup 
$\Gamma_1\times\Gamma_2$ of finite index 
such that $\Gamma_1$ and $\Gamma_2$ are infinite and 
if $X$ splits as a product $X_1\times X_2$, 
then the boundary $\partial X$ is maybe non-minimal.
The author does not have a proof of the statement: 
If $X$ splits as a product $X_1\times X_2$ 
then the boundary $\partial X$ is not minimal.
This statement seems to be true.
\end{Remark}

Here we show that 
if the answer of Question~\ref{Q1} is yes, 
then the answer of Question~\ref{Q:R} is also yes.
To prove this, we show that Question~\ref{Q:R} 
is equivalent to the following question.

\begin{Question}\label{Q:R2}
Suppose that a group $\Gamma$ acts geometrically on a CAT(0) space $X$ and 
$\Gamma$ does not contain a subgroup $\Gamma_1\times\Gamma_2$ of finite index 
where $\Gamma_1$ and $\Gamma_2$ are infinite.
Is it the case that 
the set $\{\gamma^\infty\,|\,\gamma\in\Gamma,\ o(\gamma)=\infty\}$ 
is dense in the boundary $\partial X$?
\end{Question}

\begin{Theorem}\label{thm6-2}
Questions~\ref{Q:R} and \ref{Q:R2} are equivalent.
\end{Theorem}

\begin{proof}
It is obvious that Question~\ref{Q:R} contains Question~\ref{Q:R2}.
We show that if the answer of Question~\ref{Q:R2} is yes, 
then the answer of Question~\ref{Q:R} is also yes.

Suppose that a group $\Gamma$ acts geometrically on a CAT(0) space $X$.
Let $\Gamma_1\times\dots\times\Gamma_n$ be a subgroup of finite index in $\Gamma$ 
such that each $\Gamma_i$ is infinite and each $\Gamma_i$ 
does not contain a subgroup $\Gamma_{i1}\times\Gamma_{i2}$ of finite index 
such that $\Gamma_{i1}$ and $\Gamma_{i2}$ are infinite.
Here we note that each $\Gamma_i$ is either isomorphic to $\Z$ or has finite center 
by \cite[Lemma~2.1]{Ho000}.
Hence we can suppose that 
for some number $k$, 
$\Gamma_i$ is isomorphic to $\Z$ for each $i\le k$ 
and $\Gamma_i$ has finite center for each $i>k$.

We prove by the induction on $n$.

In the case $n=1$, it is obvious.

We consider in the case $n=2$.
Then 
$\Gamma_1\times\Gamma_2$ is a subgroup of finite index in $\Gamma$ and 
$\Gamma_1\times\Gamma_2$ acts geometrically on $X$.
By Proposition~\ref{thm:splitting}, 
$X$ contains a quasi-dense subspace $X_1\times X_2$ 
such that $X_1=C(\Gamma_1x_0)$ for some $x_0\in X$ and 
$\Gamma_2$ acts geometrically on $X_2$ by projection.
Let $\alpha\in\partial X$.
Here 
$$\partial X=\partial X_1*\partial X_2
=(\partial X_1\times\partial X_2\times [-\pi,\pi])/\sim.$$
Hence $\alpha=[\alpha_1,\alpha_2,\theta]$ 
for some $\alpha_1\in\partial X_1$, 
$\alpha_2\in\partial X_2$ and $\theta\in [-\pi,\pi]$.
Now 
$\{\gamma^\infty\,|\,\gamma\in\Gamma_1,\ o(\gamma)=\infty\}$ 
is dense in $\partial X_1$ and 
$\{\delta^\infty\,|\,\delta\in\Gamma_2,\ o(\delta)=\infty\}$ 
is dense in $\partial X_2$.
Hence there exist sequences 
$\{\gamma_i\}\subset\Gamma_1$ and $\{\delta_i\}\subset\Gamma_2$ 
such that 
$\{\gamma_i^\infty\}$ converges to $\alpha_1$ and 
$\{\delta_i^\infty\}$ converges to $\alpha_2$.
Since $\langle \gamma_i,\delta_i \rangle$ is isomorphic to $\Z\times\Z$, 
by the Flat Torus Theorem (\cite[Theorem~II.7.1]{BH}), 
$\langle \gamma_i,\delta_i \rangle$ acts geometrically on 
some convex hull $C(\langle \gamma_i,\delta_i \rangle x_i)$ 
which is isometric to the Euclidean plane.
Here $C(\langle \gamma_i,\delta_i \rangle x_i)\subset X_1\times X_2$ and 
$$\{\gamma_i^{-\infty},\gamma_i^\infty,\delta_i^{-\infty},\delta_i^\infty\}
\subset \partial(C(\langle \gamma_i,\delta_i \rangle x_i)).$$
Then there exists a sequence 
$\{a_{ij}\}\subset \langle \gamma_i,\delta_i \rangle$ 
such that $\{a_{ij}^\infty\}_j$ converges to 
$[\gamma_i^\infty,\delta_i^\infty,\theta]$.
Here the sequence $\{[\gamma_i^\infty,\delta_i^\infty,\theta]\}_i$ 
converges to $\alpha$.
Hence 
$$\alpha\in\overline{\{a_{ij}^\infty\,|\,i,j\in\N \}}\subset 
\overline{\{\gamma^\infty\,|\,\gamma\in\Gamma,\ o(\gamma)=\infty\}}.$$
Thus 
$\{\gamma^\infty\,|\,\gamma\in\Gamma,\ o(\gamma)=\infty\}$ 
is dense in $\partial X$.

We consider in the case $n>2$.
Then $\Gamma_1\times\dots\times\Gamma_{n-1}\times\Gamma_n$ 
is a subgroup of finite index in $\Gamma$.
Let $\tilde{\Gamma}_1=\Gamma_1\times\dots\times\Gamma_{n-1}$ 
and $\tilde{\Gamma}_2=\Gamma_n$.
Here we can suppose that $\tilde{\Gamma}_2$ has finite center 
or each $\Gamma_i$ is isomorphic to $\Z$ for $i=1,\dots,n$.
By the inductive hypothesis and the same argument as the proof in the case $n=2$, 
we obtain that 
$\{\gamma^\infty\,|\,\gamma\in\Gamma,\ o(\gamma)=\infty\}$ is 
dense in $\partial X$.
\end{proof}

By Proposition~\ref{prop}, 
we see that Question~\ref{Q1} contains Question~\ref{Q:R2}.
Hence we obtain the following theorem from Theorem~\ref{thm6-2}.

\begin{Theorem}\label{thm6-3}
Question~\ref{Q1} contains Question~\ref{Q:R}.
\end{Theorem}

By Theorem~\ref{thm}, 
the answer of Question~\ref{Q1} is yes 
for right-angled Coxeter groups and their boundaries.
Hence, by Theorem~\ref{thm6-3}, 
we obtain that the answer of Question~\ref{Q:R} is also yes 
in this case (Corollary~\ref{cor5-1}).

%%%%%%%%%%%%%%%%%%%%%%%%%%%%%%%%%%%%%
%             REFERENCES            %
%%%%%%%%%%%%%%%%%%%%%%%%%%%%%%%%%%%%%
%

%
\end{document}